\documentclass[11pt,draft]{article}
\usepackage{mathrsfs}
\usepackage{mathrsfs}
\usepackage{}
\usepackage{amsfonts}
\usepackage{amsmath,amsfonts,mathrsfs,amssymb,color}
\usepackage{indentfirst}
\numberwithin{equation}{section}

\newtheorem{theo.}{\quad\, Theorem}[section]
\newtheorem{defi.}{\quad\, Definition}[section]
\newtheorem{lemm.}{\quad\, Lemma}[section]
\newtheorem{coro.}{\quad\, Corollary}[section]

\begin{document}

\title{On multiple positive ground state solutions for a mean curvature equation in Minkowski space}
\author{Ruyun Ma$^a$,  \ \ \ \ Yanqiong Lu$^b$ \ \ \ \ Tianlan Chen$^{c}$,\\
\small{$^{a,b,c}$Department of Mathematics, Northwest Normal University,
Lanzhou\ 730070, P R China}}
\date{} \maketitle
\footnote[0]{E-mail: mary@nwnu.edu.cn (R. Ma), linmu8610@163.com
(Y. Lu), chentianlan511@126.com(T. Chen).
} \footnote[0] {$^*$Supported by the
NSFC (No.11361054), SRFDP(No.20126203110004) and Gansu provincial National Science Foundation of China (No.1208RJZA258).}
 \begin{abstract}

 In this paper, we show how changes in the sign of nonlinearity leads to multiple radial ground state solutions of the
mean curvature equation
$
\nabla\cdot \Big[\frac{\nabla u}{\sqrt{1-|\nabla u|^2}}\Big] +\lambda f(u)=0\ \ \text{in} \ \mathbb{R}^N
$
for sufficiently large $\lambda$ with $N\geq 2$.

\end{abstract}


{\small\bf MR(2010)\ \ \ 39A12, \ 34B18}

\baselineskip 22pt

\section{Introduction}

 Hypersurfaces of prescribed mean curvature in Minkowski space are of interest in differential geometry and in general relativity. In this paper, we are concerned with the existence and multiplicity of
such a kind of hypersurfaces which are graphs of the solution of the following problem
$$\aligned
&\nabla\cdot \Big[\frac{\nabla u}{\sqrt{1-|\nabla u|^2}}\Big] +\lambda f(u)=0,\ \ \ \  \text{in} \ \mathbb{R}^N,\\
&\ u(x)>0,\ \qquad \qquad\qquad \qquad\qquad \  \text{in} \ \mathbb{R}^N,\\
&\ u(x)\to0,\ \qquad \qquad\qquad \qquad \qquad  \   \text{as} \ |x|\to\infty,
\endaligned\eqno (1)
$$
where $f:\mathbb{R}\to\mathbb{R}$ is a local Lipschitz function with $f(0)=0$, $\lambda>0$ is a parameter and $N\geq 2$.

The differential operator we are considering has been deeply studied in the recent years, in nonlinear equations on bounded domains with various type of boundary conditions (see [1-5] and the references within) and in the whole $\mathbb{R}^N$ (see [6,7]).

 The radial solutions which only depend on $r=|x|$ of (1) satisfy the following ODE
$$\aligned
&\Big(\frac{u'}{\sqrt{1-(u')^2}}\Big)'+\frac{N-1}{r}\frac{u'}{\sqrt{1-u'^2}} +\lambda f(u)=0,\\
&u(0)=\zeta,\ \ u'(0)=0,\\
\endaligned\eqno (2)
$$
where $u\in C^2([0,+\infty])$ is now a function of $r=|x|$ alone,  and $\zeta$ has to be determined in order to have
$$\lim_{r\to \infty} u(r)=0. \eqno(3)$$
The existence of the positive solution of (1) can be
interpreted in this context as the existence of a ground state solution.

Recently, Azzollini [7] proves the existence of a ground state solutions of (1) with $\lambda=1$ by the shooting method under the  assumptions:

(f1) $f(0)=0$,

(f2) $f:[0, +\infty)\to \mathbb{R}$ is locally Lipschitz,

(f3) $\exists\ \alpha:=\inf\{\xi>0|\, f(\xi)\geq 0\}>0$,

(f4) (if $N\geq 3$), $\lim\limits_{s\to\alpha^+}\frac{f(s)}{s-\alpha}>0$,

(f5) $\exists\ \gamma>0$ such that $F(\gamma):=\int^\gamma_0 f(s)ds>0$,

(f6) $f(\xi)>0$ in $(\alpha, \xi_0]$, where
$\xi_0:=\inf\{\xi\in (0, \infty)\,|\,F(\xi)>0\}.$

\noindent He proved the following

\noindent{\bf Theorem A.} ([7, Theorem 0.1]) \ If

$\bullet$ $N\geq 3$ and $f$ satisfies (f1)-(f6),

$\bullet$ $N=2$ and $f$ satisfies (f1)-(f3), (f5) and (f6),

\noindent then (1) has a radially decreasing solution with $\lambda=1$.

The shooting argument has been used in the past to find ground state
solutions to various types of equations.  For examples,  Berestycki, Lions and Peletier [8] study the existence of a ground state
solution of the Laplace equation
$$\Delta u + f(u) = 0\ \ \ \text{in}\ \mathbb{R}^N \eqno(4)$$
with $N\geq 2$. And the case $N=1$,  Berestycki and Lions [9] find the sufficient and necessary condition for the existence of the unique solution of the problem (4).
Peletier and Serrin [10] are concerned with the existence of a ground state solution of the following prescribed mean curvature equation
$$\nabla\cdot\Big[\frac{\nabla u}{\sqrt{1+|\nabla u|^2}}\Big] -\lambda u+u^q=0.$$
The shooting method consists in studying the profile of the solution of (2) as the initial value $\zeta$ varies into an interval. The main ideas is to exclude the cases in which for a finite $R>0$ either $u$ or $u'$ vanishes.

On the other hand, D\'{a}vila  del Pino and Guerra [11] find the problem
$$\Delta u-u+u^p+\lambda u^q=0\ \ \ \text{in}\ \mathbb{R}^N $$
has at least three positive decaying radial solutions if $N=3, 1<q<3$, $q<p<5$ is taken sufficiently close to $5$ and $\lambda$ is fixed sufficiently large.

Naturally, what is really interesting is to find the conditions which permit to multiple ground state solutions of (1).
Motivated above papers [6-11], this paper devotes to studying how changes in the sign of $f(s)$ leads to multiple positive radial solutions of (1).

We make the following assumptions:

(A1) $f:[0, +\infty)\to \mathbb{R}$ is locally Lipschitz with $f(0)=0$;

(A2) there exists $2n$ real numbers $0=:\beta_0<\alpha_1<\beta_1<\cdots<\alpha_n<\beta_n<\infty$
such that for $i\in \{1, \cdots, n\}$,
$$f(s)<0,\ \ \ s\in (\beta_{i-1}, \alpha_i); \ \ \ \ \ \
  f(s)>0,\ \ \ s\in (\alpha_i,\beta_i);
  $$

$(A2)'$ there exists $2n-1$ real numbers $0<\alpha_1<\beta_1<\cdots<\alpha_n<\infty$
such that for $i\in \{1, \cdots, n-1\}$,
$$f(s)<0,\ \ \ s\in (\beta_{i-1}, \alpha_i); \ \ \ \ \ \
  f(s)>0,\ \ \ s\in (\alpha_i,\beta_i); \ \ \ \ \ \ f(s)>0,\ \ \ s\in (\alpha_n,\infty);
  $$

(A3) let $F(u):=\int^u_0 f(s)ds$. Then for each $i\in\{1,2,\cdots n\}$, there exists $\xi_i\in(\alpha_i, \beta_i)$ such that $ F(\xi_i)=0$;

(A4)  for each $i\in\{1,2,\cdots n\}$, $F(\beta_{i-1})<F(\beta_{i})$;

(A5) (if $N\geq 3$) $\lim\limits_{s\to\alpha_i^+}\frac{f(s)}{s-\alpha_i}>0,\ i=1,2,\cdots, n$ and $\lim\limits_{s\to\beta_j^+}\frac{f(s)}{\beta_j-s}>0,\ j=1,2,\cdots, n-1$.

In the sequel, we will suppose that $f$ is extended in
$\mathbb{R}$ by setting $f(s)=0$ if $s\leq0$. Clearly, $f$ is locally Lipschitz continuous on $\mathbb{R}$.  The main result of the paper is the following:

\noindent{\bf Theorem 1.} \ \ If

$\bullet$ $N\geq 3$ and $f$ satisfies (A1)-(A5),

$\bullet$ $N=2$ and $f$ satisfies (A1)-(A4),

\noindent then (1) has $n$ distinct radially decreasing solutions for $\lambda>0$ is sufficiently large .

\noindent{\bf Remark 1.} \ Note that Theorem 1 is Theorem A in the case $n=1$ and $\lambda=1$.

\noindent{\bf Remark 2.} \ If we replace (A2) with (A2)$'$, then the result of Theorem 1 is also true by a similar argument with obvious changes.

\noindent{\bf Remark 3.} \ We exhibit some examples of functions $f$ satisfying our assumptions:
consider the function
 $$f(s)=(s^2-19s+18)(12s^2-s^3-27s).$$
By a simple computation, we can get that $f$ satisfies (A1)-(A5) with $\alpha_1=1, \alpha_2=9, \beta_1=3, \beta_2=18$.
From Theorem 1, there exist numbers $\zeta_i\in (\xi_i, \beta_i)$, $i=1,2$ such that for sufficiently large $\lambda$, the problem (1) has two distinct positive, decaying radial solutions.

\section{Proof of the main result}

Since we are interested in the multiplicity of ground state solutions of (1), we aim to find $n$ distinct numbers $\zeta_i\in (\xi_i, \beta_i)$, $i=1,2,\cdots, n$ such that for $\lambda>0$ is sufficiently large , the solution  $u_i\in C^2(\mathbb{R}_+)$ of the IVP:
$$\aligned
&\Big(\frac{u'}{\sqrt{1-(u')^2}}\Big)'+\frac{N-1}{r}\frac{u'}{\sqrt{1-u'^2}} +\lambda f(u)=0,\\
&u(0)=\zeta_i,\ \ u'(0)=0\\
\endaligned\eqno (5)_i
$$
has the properties:
$u_i(r)>0$ for $r\in [0, \infty)$, $u_i'(r)<0$ for $r\in (0, \infty)$ and
$$\underset{r\to+\infty}\lim u_i(r)=0.
\eqno (6)
$$

Observe that the solution of $(5)_i$ satisfies the equation
$$(r^{N-1}\phi'(u'))'=-r^{N-1}\lambda f(u),  \eqno(7)$$
where $\phi(s):=1-\sqrt{1-s^2}$ for $s\in[-1,1]$. It is easy to verify that $\phi':(-1,1)\to \mathbb{R}$ is an increasing diffeomorphism. Set $\delta>0$ and denote by $C:=C([0,\infty),\mathbb{R})$ and by $C_\delta:=C([0,\delta],\mathbb{R})$. Define the following operators
$$S:C\to C, \ \ \ \ Su(r):=\left\{
    \begin{array}{ll}
      -\frac{1}{r^{N-1}}\int_0^rt^{N-1}u(t)dt,& \text{if}\ r>0, \\
      0,\ \  & \text{if}\ r=0, \\
    \end{array}
  \right.
 $$
and $K:C\to C, \ K(u)(r)=\int_0^ru(t)dt.$

For every $\zeta_i\in\mathbb{R}$, define the translation operator $T_{\zeta_i}: C\to C$ such that $T_{\zeta_i}(u)=\zeta_i+u$. Moreover, consider the Nemytskii operators associated to $f$ and $(\phi')^{-1}$,
$$N_f:C\to C,\ \ \ \ N_f(u)(r)=f(u(r)),$$
$$N_{(\phi')^{-1}}:C\to C,\ \ \ \ N_{(\phi')^{-1}}(u)(r)=(\phi')^{-1}(u(r)).$$
Set $\rho_i>0$ and denote with $B_{\rho_i}:=\{u\in C_\delta\,|\,\|u\|_\infty\leq \rho_i\}$. We set the following fixed point problem: for any $\zeta_i\in\mathbb{R}$ we want to find $u\in\zeta_i+B_{\rho_i}$ such that
$$u=T_{\zeta_i}\circ K\circ N_{(\phi')^{-1}}\circ S\circ (\lambda N_f(u)). \eqno(8)$$
Since $(\phi')^{-1}$ and $f$ are respectively Lipschitz and locally Lipschitz, Banach-Caccioppoli fixed point theorem guarantees the existence of a sufficiently small $\delta>0$ such that the function $u_i(\lambda):=u(\zeta_i, r)\in \zeta_i+B_{\rho_i}$ is a solution of (8). It is easy to see that $u_i$ is a local solution of the Cauchy problem $(5)_i$.

Let $R_{\zeta_k}>0$ be such that $[0,R_{\zeta_k})$ is the maximal interval where the function $u_k$ is defined, here $k=1,2,\cdots, n$.
Multiplying $(5)_k$ by $u'_k$ and integrating over $(0,r)$ we obtain the following equality for any $r\in(0,R_k)$:
$$ H(u'_k(r))+(N-1)\int_0^r\frac{[u'_k( s)]^2}{s\sqrt{1-[u'_k( s)]^2}}=\lambda[F(\zeta_k)-F(u_k(r))], \eqno(9)$$
where $H(t):=\frac{1-\sqrt{1-t^2}}{\sqrt{1-t^2}}$.

For each $k\in\{1,2,\cdots,n\}$, let $I_k=(\alpha_k,\beta_k)$, and take $\zeta_k\in I_k$. By (A2) and (A4), for every $s\in (\beta_{k-1},\beta_k]$, we have $F(s)\geq F(\alpha_k)$. Thus from (9), we deduce that $H(u'(r))$ is bounded as far as $\beta_{k-1}\leq u\leq \beta_k$. Obviously, since $f(u_k(0))=f(\zeta_k)>0$, from Eq.$(5)_k$ we deduce that $u''_k(0)<0$ and this implies that there exists $\sigma>0$ such that
$$u'_k(r)<0\ \ \ \text{and}\ \ 0<u_k(r)<\zeta_k\ \ \ \text{for}\ 0<r<\sigma.
$$
Set
$$\bar{R}_{\zeta_k}:=\left\{
    \begin{array}{ll}
      \inf\{r\in(0,R_{\zeta_k})\,|\, u'_k(r)\geq 0\},\ & \text{if}\ u'_k(r)=0\ \ \text{for some}\ r\in (0, {R}_{\zeta_k}), \\
      +\infty, \ &\text{otherwise}. \\
    \end{array}
 \right. \eqno(10)
 $$
From [7, Remark 1.1], it follows that $0<\sigma\leq \bar{R}_{\zeta_k}\leq +\infty$ and  for every $r\in(0,\bar{R}_{\zeta_k})$,
$$\exists\ \varepsilon>0\ \text{such that for any}\ r\in(0,\bar{R}_{\zeta_k}),\ |u_k'(r)|\leq 1-\varepsilon. \eqno(11)$$
In particular, $\bar{R}_{\zeta_k}=+\infty$ implies ${R}_{\zeta_k}=+\infty$.

Define the following two classes of intervals
$$I_k^+:=\big\{\zeta_k\in I_k\,\big|\, \exists\ R'_k\leq R_{\zeta_k} \ \text{such that} \; u_k(r)>0, \;u'_k(r)<0\ \text{for}\  r\in (0, R'_k), u'_k(R'_k)=0\big\},
$$
and
$$I_k^-:=\big\{\zeta_k\in I_k\,\big|\, \exists\ R'_k\leq R_{\zeta_k} \ \text{such that} \; u_k(r)>0, \;u'_k(r)<0\ \text{for}\  r\in (0, R'_k), u_k(R'_k)=0\big\}.
$$
We will prove that
 the sets $I_k^+$ and $I_k^-$ are non-empty, disjoint and open, $k=1,2,\cdots n.$ Moreover, $I_k^+$ and $I_k^-$ do not cover $I_k$.


\vskip 3mm

\noindent{\bf Lemma 2.1.} Assume that $R_{\zeta_k}=+\infty$. For any fixed $\lambda>0$ and $k\in\{1,2,\cdots,n\}$,  $\zeta_k\in (0, \infty)$ be such that $u_k(r)>0$ for all $r\geq 0$ and
 $u'_k(r)<0$ for all $r>0$. Then the number $l=\underset{r\to\infty}\lim u_k(r)$ satisfies
 $$f(l)=0.
 $$
 Furthermore, if $f$ satisfies (A2) and (A5),
then $l=0$.

\noindent{\bf Proof.} Clearly, there exists $l=\lim\limits_{r\to+\infty} u(r)\geq 0$.  By $(5)_k$ and (11), we imply that
$$\lim\limits_{r\to+\infty}\Big(\frac{u'_k(r)}{\sqrt{1-[u'_k(r)]^2}}\Big)
=-\lambda f(l). \eqno(12)$$
Suppose that $f(l)\neq 0$, say $f(l)>0$. By simple computations, together with (11) and (12), we deduce that, definitively, $u''_k(r)<-\delta<0$ for some $\delta>0$. Of course this is not possible because of (11). Therefore, $f(l)=0$.

Now, we claim that $l=0$.

To this end, we only need to prove that for $k=1$, $l\neq \alpha_1$ and for each $k\in\{2,\cdots,n\}$, $l\neq \alpha_{i},\ i=1, 2, \cdots, k,\ \l\neq \beta_{j}$, $j=1, 2, \cdots, k-1$.  We divide into three steps.

Step 1. We  show that for $k=1$, $l\neq \alpha_1$.

If $N=2$ and, by contraction, $l=\alpha_1$. Since for any $r>0$, $\alpha_1<u_1(r)<\beta_1$,  from (7) we deduce that $r\phi'(u'_1(r))$ is decreasing in $[0,+\infty)$ and then, in particular, there exist $R_1>0$ and $\delta>0$ such that for any $r>R_1$, we have $\phi'(u'_1(r))<-\frac{\delta}{r}$. By (11) we infer that, for some $M_1>0$, we have $M_1u'_1(r)\leq \phi'(u'_1(r))$ and then
$$u'_1(r)\leq -\frac{\delta}{ M_1r}\ \ \ \text{for any}\ \ r>R_1.$$
Integrating in $(R_1,r)$ we obtain $$u_1(r)\leq u_1(R_1)-\frac{\delta}{ M_1}\ln\Big(\frac{r}{R_1}\Big)\to -\infty  \ \ \text{as}\ r\to+\infty,$$
which contradicts $l=\alpha_1$.

If $N\geq 3$, and suppose on the contrary
that $l=\alpha_1$, then computing in $(5)_1$, we have that the following equality holds in $(0,+\infty)$:
$$\frac{u''_1}{[1-(u'_1)^2]^{\frac{3}{2}}}=-
\frac{N-1}{r}\frac{u'_1}{\sqrt{1-(u'_1)^2}}-\lambda f(u_1).$$
Taking into account (11), there exists $\delta>0$ such that $\delta\leq \sqrt{1-(u'_1)^2}\leq 1$. We deduce that
$$u''_1=-\frac{N-1}{r}u'_1[1-(u'_1)^2]-\lambda f(u_1)[1-(u'_1)^2]^{\frac{3}{2}}\leq
-\frac{N-1}{r}u'_1-\delta^3 \lambda f(u_1), \eqno(13)$$
where we have used the fact that $u'_1<0$ and $f(u_1)>0$. Now we proceed as in [7,8], repeating the arguments for completeness. If we set $v=r^{\frac{N-1}{2}}(u_1-\alpha_1)$, by (13) we get the following estimate
 $$
v''\leq\Big\{\frac{(N-1)(N-3)}{4r^2}-\delta^3\lambda \frac{f(u_1)}{u_1-\alpha_1}\Big\}v\\
\eqno(14)
 $$
from which, in view of (A5), we deduce that $v''$ is definitively negative. Now, since $v'$ is definitively decreasing, certainly there exists $L=\lim\limits_{r\to+\infty} v'(r)<+\infty.$

However, $L$ cannot be negative, since otherwise $\lim\limits_{r\to+\infty} v(r)=-\infty$, this is a contradiction.
On the other hand, if $L\geq 0$, then we deduce that $v$ is definitively increasing and then there exists $R_1>0$ such that for any $r>R_1$, we have $v(r)>v(R_1)$. From (14) we infer that, for some positive constant $C$, $v''(r)\leq -C<0$ definitively and this implies $L=\lim\limits_{r\to+\infty}v'(r)=-\infty$, again a contradiction.

Step 2. we show that for $k=2$, $l\neq\alpha_1, \ \alpha_2$ and $l\neq \beta_1$.

By a similar argument as step 1 with $u_2(r)$ instead of $u_1(r)$,
and $v_{2,i}=r^{\frac{N-1}{2}}[u_2(r)-\alpha_i]$, $i=1,\ 2$ instead of $v$, we can deduce that $l\neq \alpha_i,\ i=1,2$. Notice that when $N=2$, we prove $l\neq \alpha_1$, by contradiction, suppose  $\lim\limits_{r\to+\infty} u_2(r)=\alpha_1$, this implies that there exists $R_2>0$ large enough such that $\alpha_1\leq u_2(r)\leq \beta_1$ for $r>R_2$, by a same argument as step 1, which deduce a contradiction.
So, we only need to show $l\neq \beta_1$.

Suppose on the contrary that $l=\beta_1$. If $N=2$, then  it follows from $\lim\limits_{r\to+\infty} u_2(r)=\beta_1$ that there exists $\tilde{R}_2>0$ large enough such that for any $r>\tilde{R}_2$, $\beta_1<u_2(r)<\alpha_2$.  If $N=2$, and $l=\beta_1$,  from (7) we deduce that $r\phi'(u'(r))$ is increasing in $[\tilde{R}_2,+\infty)$ and then, in particular, there exist $R_2>\tilde{R}_2$ and $\delta_1>0$ such that for any $r>R_2$, we have $\phi'(u'(r))>\frac{\delta_1}{r}$. By (11) we infer that, for some $M_2>0$, we have $M_2u'(r)\geq \phi'(u'(r))$ and then
$$u'(r)\geq \frac{\delta_1}{ M_2r}\ \ \ \text{for any}\ \ r>R_2.$$
Integrating in $(R_2,r)$ we obtain $$u(r)\geq u(R_2)+\frac{\delta_1}{ M_2}\ln(\frac{r}{R_2})\to +\infty  \ \ \text{as}\ r\to+\infty,$$
which contradicts $l=\beta_1$.

If $N\geq 3$, then computing in $(5)_2$, we have that the following equality holds in $(0,+\infty)$:
$$\frac{u''_2}{[1-(u'_2)^2]^{\frac{3}{2}}}=-
\frac{N-1}{r}\frac{u'_2}{\sqrt{1-(u'_2)^2}}-\lambda f(u_2).$$
Taking into account (11), there exists $\delta_2>0$ such that $\delta_2\leq \sqrt{1-(u'_2)^2}\leq 1$. We deduce that
$$u''_2=-\frac{N-1}{r}u'_2[1-(u'_2)^2]-\lambda f(u_2)[1-(u'_2)^2]^{\frac{3}{2}}\geq
-\delta_2^2\frac{N-1}{r}u'_2-\delta_2^3 \lambda f(u_2), \eqno(15)$$
where  we have used the fact that $u'_2<0$ and $f(u_2)<0$ on $[\tilde{R}_2,\infty)$. 
 
 \textcolor{blue}{Let\,$w(r)=r^{\frac{\delta^2(N-1)}{2}}(\beta_1-u(r))$. It follows that
$$w'(r)=\frac{\delta^2(N-1)}{2}r^{\frac{\delta^2(N-1)-2}{2}}(\beta_1-u(r))
-r^{\frac{\delta^2(N-1)}{2}}u'(r),$$ $$w''(r)=-r^{\frac{\delta^2(N-1)}{2}}u''(r)-\delta^2(N-1)r^{\frac{\delta^2(N-1)-2}{2}}u'(r)
+\frac{\delta^2(N-1)[\delta^2(N-1)-2]}{4r^2}r^{\frac{\delta^2(N-1)}{2}}(\beta_1-u(r)).$$
This together with inequality (15) implies that
$$
\aligned
w''(r)\leq& r^{\frac{\delta^2(N-1)}{2}}[\delta^2\frac{N-1}{r}u'(r)+\delta^3\lambda f(u)]-\delta^2(N-1)r^{\frac{\delta^2(N-1)-2}{2}}u'(r)\\
&
+\frac{\delta^2(N-1)[\delta^2(N-1)-2]}{4r^2}r^{\frac{\delta^2(N-1)}{2}}(\beta_1-u(r))\\
=&\delta^2\frac{N-1}{r}r^{\frac{\delta^2(N-1)}{2}}u'(r)
-\delta^2\frac{N-1}{r}r^{\frac{\delta^2(N-1)}{2}}u'(r)\\
&+[\delta^3\lambda \frac{f(u)}{\beta_1-u}+\frac{\delta^2(N-1)[\delta^2(N-1)-2]}{4r^2}]w(r)\\
=&[\delta^3\lambda \frac{f(u)}{\beta_1-u}+\frac{\delta^2(N-1)[\delta^2(N-1)-2]}{4r^2}]w(r),
\endaligned \eqno(16)$$
from which, in view of (A5), we deduce that $w''$ is definitively negative. Now, since $w'$ is definitively decreasing, certainly there exists $L=\lim\limits_{r\to+\infty} w'(r)<+\infty.$
However, $L$ cannot be positive, since otherwise $\lim\limits_{r\to+\infty} w(r)=+\infty$, this is a contradiction.
On the other hand,
 if $L\leq0$, then
 we deduce that there exists $R_2>0$ such that
 $$w'(r)\leq 0, \ \ \ \ \  r>R_2,
 $$ and $w$ is definitively decreasing.
 Hence, there exist two constants $R_*$ and $\sigma$ with $R_\ast>R_2$ and $\sigma>0$, such that
 $$w(r)<0, \   w''(r)\leq \sigma^2 w(r),\ \ \ \ \ \ r\in[R_\ast,+\infty).
 $$
\
 Set $b_1:=w(R_\ast),  b_2:=w'(R_\ast)$. Then  $b_1<0, b_2\leq 0$.
Let us consider the initial value problem
$$x''(r)=\sigma^2 x(r),\ r\in(R_\ast,\infty),\ \ \ x(R_\ast)=b_1,\ \ x'(R_\ast)=b_2.$$
 Its unique solution can be explicitly given by
$$x(r)=\frac{b_1-b_2}{2} e^{-\sigma(r-R_\ast)}+\frac{b_1+b_2}{2} e^{\sigma(r-R_\ast)},\ \ r\in (R_\ast,\infty).$$
Let $z(r)=x(r)-w(r)$. Then
$$z''(r)\geq \sigma^2z(r),\ \ r\in (R_\ast,\infty),\ \ \ z(R_\ast)=0,\ \ z'(R_\ast)=0. $$
Let
$$M(r):=z''(r)-\sigma^2 z(r), \ \ \ \  r\in (R_\ast,\infty).
$$
Then 
$$z''(r)-\sigma^2 z(r)=M(r),\ \ r\in (R_\ast,\infty),\ \ \ z(R_\ast)=0,\ \ z'(R_\ast)=0,
$$
which has a unique solution
$$z(r)=\frac 1{2\sigma} \int^r_{R_\ast} [e^{\sigma(r-s)}-e^{-\sigma(r-s)}]M(s)ds,\ \ \ \ \ \  r\in (R_\ast,\infty).
$$
Obviously, $z(r)\geq 0$ for $r\in (R_\ast,\infty)$, which implies
$$x(r)\geq w(r),\ \ \ \ \ r\in (R_\ast,\infty),$$
i. e.
$$\frac{b_1-b_2}{2} e^{-\sigma(r-R_\ast)}+\frac{b_1+b_2}{2} e^{\sigma(r-R_\ast)}\geq r^{\frac{\delta^2(N-1)}{2}}(\beta_1-u(r)),\ \ \ \  r\in (R_\ast,\infty).$$
However, this is impossible since
$$\frac{b_1+b_2}{2}<0, \ \ \beta_1-\zeta_2<\beta_1-u(r)\leq 0, \ \ \ \lim_{r\to \infty}\frac {e^{\sigma(r-R_\ast)}}{r^{\frac{\delta^2(N-1)}{2}}}=+\infty.
$$
}
Therefore, $l\neq \beta_1$.

\vskip 3mm
Step 3. We claim that for each $k\in\{3, 4, \cdots, n\}$, $l\neq\alpha_i, \ i=1, 2, \cdots, k$ and $l\neq \beta_j$, $j=1, 2, \cdots, k-1$.

Toward this end, we only need to repeat the arguments of step 1 and step 2 with $u_k(r)$ instead of $u_1(r)$,
and  $v_{k,i}=r^{\frac{N-1}{2}}[u_k(r)-\alpha_i]$, $i=1,\ 2,\cdots,k$ instead of $v$, and the proof of $l\neq \beta_1$ with
$u_k$ instead of $u_2$ and $w_{k,j}(r)=r^{\frac{N-1}{2}}[\beta_j-u_k(r)],\ j=1,2,\cdots,k-1$  instead of $w$.
\hfill$\Box$

\vskip 2mm

\noindent{\bf Lemma 2.2.} For any fixed $\lambda>0$ and let $k\in\{1,2,\cdots,n\}$, $I_k^+\neq\emptyset$.

\noindent{\bf Proof.}
 Let $\zeta_k\in (\alpha_k, \xi_k]$. By (A2) and (A3), $F(\zeta_k)<0$. Because of (9) and the definition of $\xi_k$, it is clear to see that $F(u(r))<F(\zeta_k)<0$ for any $r\in (0,R_{\zeta_k})$. As a consequence, by the fact $f(\zeta_k)>0$ in $(\alpha_k,\xi_k]$ we have that there exists $m_k>0$ such that $$0<m_k<u_k(r)<\zeta_k. \eqno(17)$$
Suppose on the contrary that $\zeta_k\not\in I_k^+$, then $\bar{R}_{\zeta_k}=+\infty$ implies $R_{\zeta_k}=+\infty$. So $u'(r)<0$ for any $r>0$, by Lemma 2.1 we get a contradiction with (17).
   \hfill{$\Box$}

\vskip 2mm

Next, we will prove that $I_k^-$ is not empty, we need some preliminary  results.

For each $i\in\{1,2,\cdots,n\}$, consider the problem
$$\aligned
&\nabla\cdot\Big[\frac{\nabla u}{\sqrt{1-|\nabla u|^2}}\Big] +\lambda f(u)=0,\ \ \ \  \text{in} \ B_{\rho},\\
&\ u=0,\ \qquad \qquad \qquad \qquad \qquad \quad \ \ \text{on} \ \partial B_{\rho}.
\endaligned\eqno (18)
$$
Recall the definition of $\beta_i$, we replace $f$ in (18) by
$$f_i(s)=\left\{
           \begin{array}{ll}
            f(s), & \text{if}\ s\leq \beta_i, \\
            f(\beta_i), & \text{if}\  s>\beta_i. \\
           \end{array}
         \right. \eqno(19)
$$
As in [3, 7], we use a variational approach to (18).

Set $W_{\rho}:=W^{1,\infty}((0,\rho),\mathbb{R})$. It is well known that $W_{\rho}\hookrightarrow C_{\rho}$. Define
$$K:=\{u\in W_{\rho}\,|\, \|u'\|_\infty\leq 1,\ \ u(\rho)=0\}$$
and $$\Psi(u):=
\left\{
    \begin{array}{ll}
      \int_0^{\rho}r^{N-1}\Big(1-\sqrt{1-(u')^2}\Big)dr, &\text{if}\ u\in K, \\
      +\infty,\   &\text{if}\ u\in W_{\rho}\backslash K.\   \\
    \end{array}
  \right.
 $$
 For any $u\in W_{\rho}$, we set $$J_i(\lambda,u):=\Psi(u)-\lambda\int_0^{\rho} r^{N-1} F_i(u)dr.$$
It is easy to verify that the functional $J_i(\lambda, \cdot)$ is a Szulkin's functional (see [12]) so that, by [12, Proposition 1.1], we have that if $u\in W_{\rho}$ is a local minimum of $J_i(\lambda, \cdot)$, then it is a Szulkin critical point and for any $v\in K$ it solves the inequality
$$\int_0^{\rho}r^{N-1}(\phi(v')-\phi(u'))dr-\lambda\int_0^{\rho}r^{N-1}f_i(u)(v-u)dr\geq 0, \eqno (20)$$
where we recall that $\phi$ is defined in (7). By a similar argument from [7, 13], we obtain the following lemma.

\noindent{\bf Lemma 2.3.}\ For all $\lambda>0$, if $v_i(\lambda,\cdot)\in K$ is a local minimum for $J_i(\lambda, \cdot)$, then $v_i(\lambda,|x|)$ is a classical solution of (18) for each $i\in\{1,2,\cdots, n\}$.

\noindent{\bf Lemma 2.4.}\ For all $\lambda>0$, $\forall\ \rho>0$ and for each $i\in\{1,2,\cdots, n\}$, there exists $v_i(\lambda,\cdot)\in K$ such that $J_i(\lambda, \cdot)$ attains its local minimum at $v_i(\lambda,\cdot)$.

Moreover, $v_i(\lambda,\cdot)$ is a classical nontrivial solution of (18) and satisfies $0\leq v_i(\lambda,\cdot)\leq \beta_i$.

\noindent{\bf Proof.}\ As a first step, we show that $J_i(\lambda,\cdot)$ is bounded below and achieves its infimum.

Observe that $\forall\ v\in K, \|v\|_\infty\leq \rho$. As a consequence, it is easy to see that $J_i(\lambda,\cdot)$ is bounded below. Consider $\{v_{i,k}\}_{k=1}^{\infty}\in W_{\rho}$ a minimizing sequence. Of course we can assume $v_{i,k}\in K$ for any $k\geq 1$. By the Ascoli Arzel\`{a}
theorem, there exists a subsequence, relabeled $\{v_{i,k}\}_{k=1}^{\infty}$, and a continuous function $v_{i}^\ast$ such that $$v_{i,k}\to v_{i}^\ast\ \ \ \text{uniformly in}\ \ [0,\rho]. \eqno (21)$$
To prove that $v_{i}^\ast$ is in $K$, we just observe that, for any $x, y\in [0,\rho]$ with $x\neq y$, we have
$$\lim\limits_{k}\frac{v_{i,k}(x)-v_{i,k}(y)}{x-y}
=\frac{v_{i}^\ast(x)-v_{i}^\ast(y)}{x-y},$$
and then also $v_{i}^\ast$ has Lipschitz constant $1$. By (21) and [13, Lemma 1], it deduce that $\Psi(v_{i}^\ast)\leq \liminf\limits_{k}\Psi(v_{i,k}).$
Then, again by (21), we have
$$J_i(\lambda, v_{i}^\ast)\leq c_{i,0},$$
where $c_{i,0}=\inf\limits_{v\in W_{\rho}} J_i(\lambda, v)$.

Now we claim that if $\rho>0$ is sufficiently large, then
$c_{i,0}<0.$
Consider the following function defined for $\rho>2\gamma_i$,
$$\omega_{\rho}(r)=\left\{
           \begin{array}{cc}
           \gamma_i, & \text{in}\ \ [0,\rho-2\gamma_i], \\
            \frac{\rho-r}{2}, & \text{in}\  \ [\rho-2\gamma_i,\rho]. \\
           \end{array}
         \right.$$
Of course $\omega_{\rho} \in K$. Moreover
$$\aligned
J_i(\lambda,\omega_{\rho})\leq&\frac{1}{2}\int_{\rho-2\gamma_i}^{\rho}(2-\sqrt{3})s^{N-1}ds
-F(\gamma_i)\frac{(\rho-2\gamma_i)^N}{N}+\frac{1}{N}\max\limits_{0\leq s\leq \gamma_i}|F(s)|[(\rho)^N-(\rho-2\gamma_i)^N]\\
\leq &C_1[\rho^N-(\rho-2\gamma_i)^N]-
\frac{F(\gamma_i)(\rho-2\gamma_i)^N}{N}\\
\leq&C_2\rho^{N-1}-C_3\rho^N<0\ \ \ \text{as}\ \rho>2\gamma_i\ \ \text{sufficiently large},
\endaligned
$$
where $C_1,\ C_2$ and $C_3$ are suitable positive constants. The claim is an obvious consequence of the previous chain of inequalities.
This together with Lemma 2.3 yields the conclusion.
\hfill$\Box$

\vskip 2mm

We will use a similar method in [14, Lemma 2.5] to obtain an important lemma.

\noindent{\bf Lemma 2.5.}\ If $\lambda>0$ is sufficiently large, then
$\sup\{v_{i+1}(\lambda,r)\,|\, r\in B_\rho\}>\beta_{i}$ and consequently, $v_{i+1}(\lambda,\cdot)\neq v_i(\lambda,\cdot)$, $i\in\{1,2,\cdots, n-1\}$.

\noindent{\bf Proof.}\  To this end, we only need to show that there exists $w\in K$ such that $J_{i+1}(\lambda, w)<J_{i+1}(\lambda, v)$
for all $v\in K$ satisfying $0\leq v\leq \beta_i$.

First of all, we show that for $\lambda>0$ is sufficiently large, then
$\sup\{v_{2}(\lambda,r)\,|\, r\in B_\rho\}>\beta_{1}$ and subsequently $v_{2}(\lambda,\cdot)\neq v_1(\lambda,\cdot)$.

 Let $\varrho_0=\inf\{F(\beta_2)-F(v(r)): r\in \bar{B}_{\rho}\ \text{and}\ 0\leq v\leq \beta_1\}$. Then $\varrho_0>0$ as $f$ satisfies the condition (A4).
If $v\in K$ satisfies $0\leq v\leq \beta_1$, then
$$\aligned
&\int_0^{\rho}r^{N-1} F_2(v(r))dr=\int_0^{\rho}r^{N-1}F(u(r))dr\\
\leq&\int_0^{\rho}r^{N-1} F(\beta_2)dr-\frac{\rho^N}{N}\varrho_0
=F(\beta_2)\frac{\rho^N}{N}-\varrho_0\frac{\rho^N}{N}.
\endaligned \eqno(22)
$$
On the other hand, let $\rho>2\beta_2$,
consider the following function
$$w_{\rho}(r)=\left\{
           \begin{array}{cc}
          \beta_2 , & \text{in}\ \ [0,\rho-2\beta_2], \\
           \frac{\rho-r}{2} , & \text{in}\  \ [\rho-2\beta_2,\rho]. \\
           \end{array}
         \right.$$
Obviously, $w_{\rho}\in K$ 
and
$$\aligned
\int_0^{\rho}r^{N-1} F(w_{\rho}(r))dr=&\int_0^{\rho-2\beta_2} r^{N-1} F(\beta_2)dr+\int_{\rho-2\beta_2}^{\rho} r^{N-1}F\big(\frac{\rho-r}{2}\big)dr\\
=&\int_0^{\rho} r^{N-1} F(\beta_2)dr-\int_{\rho-2\beta_2}^{\rho} r^{N-1} F(\beta_2)dr+\int_{\rho-2\beta_2}^{\rho} r^{N-1}F\big(\frac{\rho-r}{2}\big)dr\\
\geq &F(\beta_2)\frac{\rho^N}{N}-2\sup\limits_{u\in[0,\beta_2]}|F(u)|
\frac{\rho^N-(\rho-2\beta_2)^N}{N}
.\\
\endaligned\eqno(23)
$$
By (22) and (23) we can choose and fix $\rho>2\beta_2$ sufficiently large  so that
$$\aligned
&\int_0^{\rho}r^{N-1} F(w_{\rho}(r))dr-\int_0^{\rho}r^{N-1} F(v(r))dr\\
\geq&\int_0^{\rho}r^{N-1} F(w_{\rho}(r))dr-\int_0^{\rho}r^{N-1} F(v(r))dr\\
\geq &\varrho_0\frac{\rho^N}{N}-2\sup\limits_{v\in[0,\beta_2]}|F(v)|
\frac{\rho^N-(\rho-2\beta_2)^N}{N}\\
\geq& C_4\rho^N-C_5 \rho^{N-1}>0,  \qquad \forall\ 0\leq v\leq\beta_1,\\
\endaligned
$$
here $C_4,\,C_5$ are suitable positive constants. Thus, there exists $\sigma_1>0$ such that
$$\int_0^{\rho}r^{N-1} F(w_{\rho}(r))dr-\int_0^{\rho}r^{N-1} F(v(r))dr>\sigma_1$$
for all $0\leq v\leq \beta_1$.  Moreover, for such $\rho>2\beta_2$, it follows that
$$\aligned
&J_{2}(\lambda,w_\rho)-J_{2}(\lambda,v)\\
=&\int_0^{\rho} r^{N-1}[1-\sqrt{1-(w'_\rho(r))^2}]dr-\int_0^{\rho} r^{N-1}[1-\sqrt{1-(v'(r))^2}]dr\\
&-\lambda\int_0^{\rho} r^{N-1}[F(w_\rho(r))-F(v(r))]dr\\
\leq& \frac{1}{2}\int_{\rho-2\beta_2}^{\rho} (2-\sqrt{3})r^{N-1}dr -\lambda\int_0^{\rho} r^{N-1}[F(w_\rho(r))-F(v(r))]dr\\
\leq& \frac{2-\sqrt{3}}{2N}[\rho^N-(\rho-2\beta_2)^N ]-\lambda\sigma_1\\
\leq&0 \ \ \ \ \text{for\ $\lambda$\ sufficiently large}.\\
\endaligned
$$

Hence, for such $\lambda$, the local minimum of $J_2(\lambda,\cdot)$ cannot be attained at any $v\in W_\rho$ such that $0\leq v\leq \beta_1$. Therefore, $\sup\{v_2(\lambda,r)\,|\, r\in B_\rho\}>\beta_1$ and so $v_2(\lambda,\cdot)\neq v_1(\lambda,\cdot)$.

By the same argument with obvious changes, we can obtain that
$\sup\{v_{i+1}(\lambda,r)\,|\, r\in B_\rho\}>\beta_{i}$ and $v_{i+1}(\lambda,\cdot)\neq v_i(\lambda, \cdot)$, $i\in\{1,2,\cdots, n-1\}$ for $\lambda$ sufficiently large.
\hfill$\Box$

\vskip 2mm

From Lemma 2.3 to Lemma 2.5, it deduce that for any fixed $\rho>0$ large enough and $k\in\{2,3,\cdots,n\}$, the problem (18) with $f_k$ instead of $f$ has $k$ distinct nontrivial solutions and  the $k$-th solution $v_k$ satisfying $\sup\{v_{k}(\lambda,r)\,|\, r\in B_\rho\}>\beta_{k-1}$ with $\lambda>0$ sufficiently large.

\noindent{\bf Lemma 2.6}\  Let $\lambda>0$ be sufficiently large and $i\in\{1,2,\cdots,n\}$. Then $I_i^-\neq\emptyset$.

\noindent{\bf Proof} \ From Lemma 2.4, it follows that $J_i(\lambda, \cdot)$ is bounded below and achieves its infimum. Moreover,  if $\rho>0$ is sufficiently large, then
$c_{i,0}<0.$

Now choose $\rho_{i}>0$ large enough such that there exists $u_{i}:=v_i(\lambda,\cdot)\in K_i$ satisfying
$J_i(u_{i})=c_{i,0}<0$ and $\sup\{u_i(r)\,|\,r\in B_{\rho_i}\}>\beta_{i-1}$. Set $\tilde{\zeta}_i=u_{i}(0)$. Then the value $\tilde{\zeta}_i\in (\alpha_i,\beta_i)$. Indeed, by Lemma 2.3 and Lemma 2.5,
$u_{i}(|\cdot|)$ is a classical solution of (18) with $\rho_i$ instead of $\rho$, and then $u_{i}$ is a local solution of $(5)_i$ with $\zeta_i=\tilde{\zeta}_i$ and $f_i$ instead of $f$. If $\tilde{\zeta}_i\leq\alpha_i$, then $\tilde{\zeta}_i\in(\beta_{i-1},\alpha_i]$ such that $F(\tilde{\zeta}_i)\leq 0$ leads to an obvious contradiction to (9)
computed in $r=\rho_{i}$. On the other hand, $\tilde{\zeta}_i$ can not be greater than $\beta_i$, since in this case, by (19), the unique solution of the Cauchy problem $(5)_i$ would be the constant function $u_i(r)=\tilde{\zeta}_i$.

By contradiction, suppose that $\tilde{\zeta}_i\not\in I_i^-$. Since we can assume $u_{i}(r)>0$ in $[0,\rho_{i})$, otherwise we consider the function $u_{i}$ restricted to the interval $[0,R_i')$, where $R_i':=\inf\{r>0\,|\, u_{i}(r)=0\}$, our contradiction assumption implies that $\bar{R}_{\tilde{\zeta}_i}\in (0,\rho_{i})$ (the definition of $\bar{R}_{\tilde{\zeta}_i}$ is given in (10)).

Computing (9) for $r=\bar{R}_{\tilde{\zeta}_i}$ and for $r=\rho_{i}$, we respectively have
$$(N-1)\int_0^{\bar{R}_{\tilde{\zeta}_i}}
\frac{(u'_i(s))^2}{s\sqrt{1-(u'_i(s))^2}}ds
=\lambda[F(\tilde{\zeta}_i)-F(u_i(\bar{R}_{\tilde{\zeta}_i}))], \eqno(24)$$
$$H(u'_i(\rho_{i}))+(N-1)\int_0^{\rho_{i}}
\frac{(u'_i(s))^2}{s\sqrt{1-(u'_i(s))^2}}ds
=\lambda F(\tilde{\zeta}_i). \eqno(25)$$
Subtracting (24) from (25), we obtain
$$H(u'_i(\rho_{i}))+(N-1)\int_{\bar{R}_{\tilde{\zeta}_i}}^{\rho_{i}}
\frac{(u'_i(s))^2}{s\sqrt{1-(u'_i(s))^2}}ds
=\lambda F(u_i(\bar{R}_{\tilde{\zeta}_i})), $$
which implies that $F(u_i(\bar{R}_{\tilde{\zeta}_i}))>0$.

Since $u_i'(r)<0$ for any $r\in (0,\bar{R}_{\tilde{\zeta}_i})$, we have that $u_i''(\bar{R}_{\tilde{\zeta}_i})\geq 0$ and then from the equation of $(5)_i$, it follows that $f(u_i(\bar{R}_{\tilde{\zeta}_i}))\leq 0$. Since $f$ is positive in $I_i$ and $0<u_i(\bar{R}_{\tilde{\zeta}_i})<\tilde{\zeta}_i<\beta_i$, certainly $u_i(\bar{R}_{\tilde{\zeta}_i})\in(\beta_{i-1},\alpha_i]$. From this we deduce that $F(u_i(\bar{R}_{\tilde{\zeta}_i}))<0$ and then the contradiction is obtained.
\hfill$\Box$

\noindent{\bf Lemma 2.7}\ For any fixed $\lambda>0$ sufficiently  large and let $k\in\{1,2,\cdots,n\}$, $I_k^-$ and $I_k^+$ are open and disjoint.

\noindent{\bf Proof} \
By contradiction, suppose $\bar{\zeta}_k\in I_k^+\cap I_k^-$.  Then,
since the solution of $(5)_k$ with ${\zeta}_k=\bar{\zeta}_k$ is
such that $u_k(R'_{{\zeta}_k})=u'_k(R'_{{\zeta}_k})=0$,
by uniqueness theorem, $u=0$ is the unique solution of the Cauchy problem
$$\aligned
&\Big(\frac{u'_k}{\sqrt{1-(u'_k)^2}}\Big)'+
\frac{N-1}{r}\frac{u'_k}{\sqrt{1-(u'_k)^2}} +\lambda f(u)=0,\\
&u(R'_{{\zeta}_k})=0,\ \ \ \ u'(R'_{{\zeta}_k})=0.
\endaligned
$$
Finally,  by continuous dependence on the initial datum, it is easy to see that $I_k^+$ and $I_k^-$ are open sets.
\hfill$\Box$

By Lemma 2.2, Lemma 2.6 and Lemma 2.7,  for $\lambda>0$ is sufficiently large, we can take $\zeta_k\in I_k\backslash(I_k^+\cup I_k^-)$ such that $u_k(r)$ is defined on $[0,\infty)$ and, $u_i\neq u_j,\ \ i\neq j$. By Lemma 2.1, $\lim\limits_{r\to+\infty} u_k(r)=0$. As a consequence, the problem (1) has $n$ distinct positive, decaying radial solution $u_k$, $k=1,2,\cdots, n$.

\vskip 5mm

\centerline {\bf REFERENCES}
\begin{description}

 \item{[1]} C. Bereanu, P.Jebelean, J. Mawhin, Radial solutions for some nonlinear problems involving mean curvature operators in Euclidean and Minkowskis paces, Proc. Amer. Math. Soc. 137 (2009), 171-178.

 \item{[2]} C. Bereanu, P. Jebelean, J. Mawhin, Radial solutions for Neumann problems involving mean curvature operators in Euclidean and Minkowskis paces, Math. Nachr. 283 (2010), 379-391.

 \item{[3]} C. Bereanu, P. Jebelean, P.J. Torres, Positive radial solution for Dirichlet problems with mean curvature operators in Minkowski space, J. Funct. Anal. 264 (2013), 270-287.

 \item{[4]} C. Bereanu, P. Jebelean, P.J. Torres, Multiple positive radial solutions for a Dirichlet problem involving the mean curvature operator in Minkowski space, J. Funct. Anal. 265 (2013), 644-659.

 \item{[5]}  C. Corsato, F. Obersnel, P. Omari, S. Rivetti, Positive solutions of the Dirichlet problem for the prescribed mean curvature equation in Minkowski space, J. Math. Anal. Appl. 405 (2013), 227-239.

 \item{[6]} D. Bonheure, A. Derlet, C. De Coster, Infinitely many radial solutions of a mean curvature equation in Lorentz-Minkowski space, Rend. Istit. Mat. Univ. Trieste 44 (2012), 259-284.

 \item{[7]} A. Azzollini, Ground state solution for a problem with mean curvature operator in Minkowski space, Journal of Functional Analysis 266 (2014), 2086-2095.

\item{[8]} H. Berestycki, P.L. Lions, L.A. Peletier, An ODE approach to the existence of positive solutions for semilinear problems in $\mathbb{R}^N$, Indiana Univ. Math. J. 30 (1981), 141-157.

 \item{[9]} H. Berestycki, P.L. Lions, Nonlinear scalar field equations. I. Existence of a ground state,
    Arch. Ration. Mech. Anal. 82 (1983), 313-345.

 \item{[10]} L.A. Peletier, J. Serrin, Ground states for the prescribed mean curvature equation, Proc. Amer. Math. Soc. 100 (1987), 694-700.

 \item{[11]} J. D$\acute{¡äa}$vila, M. del Pino,  I. Guerra, Non-uniqueness of positive ground states of non-linear
Schr\"{o}dinger equations, Proc. London Math. Soc. 106(3) (2013), 318-344.

 \item{[12]} A. Szulkin, Minimax principles for lower semicontinuous functions and applications to nonlinear boundary value problems, Ann. Inst. H. Poincar\`{e} Anal. NonLin\`{e}aire 3 (1986), 77-109.

 \item{[13]} H. Brezis, J. Mawhin, Periodic solution of the forced relativistic pendulum, Differential Integral Equations 23 (2010), 801-810.

\item{[14]} K. J. Brown, H. Budin, On the existence of positive solutions for a class of semilinear elliptic boundary value problems, SIAM J. Math. Anal. 10(5) (1979), 875-883.


\end{description}

\end{document}